\documentclass[11pt]{amsart} 
\IfFileExists{pdfsync.sty}{\usepackage{pdfsync}}{}
\usepackage{array}

\usepackage{amsrefs}
\usepackage{varioref}
\usepackage[all]{xy}
\CompileMatrices
\SelectTips{cm}{}

\RequirePackage{amsmath}
\RequirePackage{amsopn}
\RequirePackage{amsfonts}

\renewcommand{\MR}[1]{\relax}

\newtheorem{thm}{Theorem}[section]

\newtheorem{lemma}[thm]{Lemma}
\newtheorem{prop}[thm]{Proposition}
\newtheorem{conjecture}[thm]{Conjecture}

\theoremstyle{definition}
\newtheorem{definition}[thm]{Definition}

\theoremstyle{remark}
\newtheorem{remark}[thm]{Remark}
\newtheorem{example}[thm]{Example}

\numberwithin{equation}{section}

\newcount\hours
\newcount\minutes       
\def\timeofday{
\hours=\time
\minutes=\hours
\divide\hours by60
\multiply\hours by60
\advance\minutes by-\hours
\divide\hours by60
\ifnum\hours>9\else0\fi\the\hours:\ifnum\minutes>9\else
0\fi\the\minutes}
\def\predate{\date{\the\day\ \ifcase\month\or
  January\or February\or March\or April\or May\or June\or July\or
        August\or September\or October\or November\or
           December\fi\ \the\year\ --- \timeofday\ --- Preliminary
                  Version}}

\def\mathcs{C^{*}}
\def\cs{\ifmmode\mathcs\else$\mathcs$\fi}

\DeclareMathSymbol{\rtimes}{\mathbin}{AMSb}{"6F}
\def\acg{{A \rtimes_\alpha G}}

\def\Z{\mathbf{Z}}
\def\K{\mathcal{K}}
\let\ker\relax
\let\sp\relax

\DeclareMathOperator{\Res}{Res}

\DeclareMathOperator{\Ind}{Ind}

\DeclareMathOperator{\Prim}{Prim}

\DeclareMathOperator{\ker}{ker}

\DeclareMathOperator{\Aut}{Aut}
\DeclareMathOperator*{\sp}{span}

\DeclareMathOperator{\id}{id}
\def\set#1{\{\,#1\,\}}
\let\tensor=\otimes
\def\restr#1{|_{{#1}}}

\def\specnp#1{{#1}^\wedge}
\newbox\hidebox
\def\spechide#1{\setbox\hidebox=\hbox{$#1$}
\hbox to\wd\hidebox{$\box\hidebox^\wedge$\hss}}
\makeatletter
\def\labelenumi{\textnormal{(\@alph\c@enumi)}}
\def\theenumi{\@alph \c@enumi}
\newcount\charno
\def\alphapart#1{\charno=96
\advance\charno by#1\char\charno}

\makeatother
\def\<{\langle}
\def\>{\rangle}
\let\ipscriptstyle=\scriptscriptstyle
\def\lipsqueeze{{\mskip -3.0mu}}
\def\ripsqueeze{{\mskip -3.0mu}}
\def\ipcomma{\nobreak\mathrel{,}\nobreak}
\newbox\ipstrutbox
\setbox\ipstrutbox=\hbox{\vrule height8.5pt
width 0pt}
\def\ipstrut{\copy\ipstrutbox}
\def\lip#1<#2,#3>{\mathopen{\relax_{\ipstrut\ipscriptstyle{
#1}}\lipsqueeze
\langle} #2\ipcomma #3 \rangle}
\def\blip#1<#2,#3>{\mathopen{\relax_{\ipstrut
\ipscriptstyle{ #1}}\lipsqueeze\bigl\langle} #2\ipcomma #3 \bigr\rangle}
\def\rip#1<#2,#3>{\langle #2\ipcomma #3
\rangle_{\ripsqueeze\ipstrut\ipscriptstyle{#1}}}
\def\brip#1<#2,#3>{\bigl\langle #2\ipcomma #3
\bigr\rangle_{\ripsqueeze\ipstrut\ipscriptstyle{#1}}}
\def\angsqueeze{\mskip -6mu}
\def\smangsqueeze{\mskip -3.7mu}
\def\trip#1<#2,#3>{\langle\smangsqueeze\langle #2\ipcomma #3
\rangle\smangsqueeze\rangle_{\ripsqueeze\ipstrut\ipscriptstyle{#1}}}
\def\btrip#1<#2,#3>{\bigl\langle\angsqueeze\bigl\langle #2\ipcomma
#3
\bigr\rangle
\angsqueeze\bigr\rangle_{\ripsqueeze\ipstrut\ipscriptstyle{#1}}}
\def\tlip#1<#2,#3>{\mathopen{\relax_{\ipstrut\ipscriptstyle{
#1}}\lipsqueeze \langle\smangsqueeze\langle} #2\ipcomma #3
\rangle\smangsqueeze\rangle}
\def\btlip#1<#2,#3>{\mathopen{\relax_{\ipstrut\ipscriptstyle{
#1}}\lipsqueeze
\bigl\langle\angsqueeze\bigl\langle} #2\ipcomma #3
\bigr\rangle\angsqueeze\bigr\rangle}

\def\ip(#1|#2){(#1\mid #2)}
\def\bip(#1|#2){\bigl(#1 \mid #2\bigr)}
\def\Bip(#1|#2){\Bigl( #1 \bigm| #2 \Bigr)}
\newcommand{\aga}{(A,G,\alpha)}
\newcommand{\acgp}{A\rtimes_{\alpha}G_{P}}
\newcommand{\ach}{A\rtimes_{\alpha}H}
\newcommand{\indgpg}{\Ind_{G_{P}}^{G}}

\newcommand{\EHI}{EHI}
\renewcommand{\H}{\mathcal{H}}

\newcommand{\X}{\mathsf{X}}
\newcommand{\lt}{\operatorname{lt}}
\newcommand{\E}{C_{0}(G/H,A)\rtimes_{\lt\tensor\alpha}G}

\def\sme{\,\mathord{\mathop{\text{--}}\nolimits_{\relax}}\,}
\def\ib{im\-prim\-i\-tiv\-ity bi\-mod\-u\-le}
\def\ibind#1{\mathop{#1\mathord{\mathop{\text{--}}}}\!\Ind\nolimits}
\def\xind{\ibind{\X}}
\newcommand{\V}{\mathcal{V}}
\newcommand{\Hind}{L^{2}_{V}(G,\beta;\H)}
\newcommand{\sHind}{\L^{2}_{V}(G,\beta;\H)}
\usepackage{mathrsfs}
\renewcommand{\L}{\mathscr{L}}

\newcommand{\IC}{\mathcal{IC}}

\def\specnp#1{{#1}^\wedge}
\newbox\hidebox
\def\spechide#1{\setbox\hidebox=\hbox{$#1$}
\hbox to\wd\hidebox{$\box\hidebox^\wedge$\hss}}
\def\rank{\operatorname{rank}}
\def\dim{\operatorname{dim}}
\def\NN{\mathbb N}
\def\GL{\operatorname{GL}}
\def\RR{\mathbb R}
\def\TT{\mathbb T}
\def\dach{\;\widehat{\ \ }}
\def\om{\omega}
\def\Om{\Omega}
\newcommand{\hA}{\hat A}

\allowdisplaybreaks[2]

\begin{document}

\title{Inducing primitive ideals}

\date{21 October 2005}

\author[Echterhoff]{Siegfried Echterhoff}
\address{Westf\"alische Wilhelms-Universit\"at M\"unster \\
Mathematisches Institut \\
Einsteinstr. 62 \\
W-48149 M\"unster\\
Germany}
\email{echters@math.uni-muenster.de}

\author[Williams]{Dana P. Williams}
\address{Department of Mathematics \\
  Dartmouth College \\
  Hanover, NH 03755-3551 \\
  USA} 
\email{dana.williams@dartmouth.edu}

\thanks{Partly supported by the Ed Shapiro fund at Dartmouth College
  and the Deutsche Forschungsgemeinschft (SFB 478)}

\begin{abstract}
  We study conditions on a $C^*$-dynamical system $\aga$ under which
  induction of primitive ideals (resp.\ irreducible representations)
  from stabilizers for the action of $G$ on the primitive ideal space
  $\Prim(A)$ give primitive ideals (resp.\ irreducible
  representations) of the crossed product $\acg$. The results build on
  earlier results of Sauvageot \cite{sau:jfa79} and
  others, and will
  correct a (possibly overly optimistic) statement of the first author in
  \cite{ech:mams96}. In an appendix, the first author takes the
  opportunity to fill a gap in the proof of another result in
  \cite{ech:mams96}.
\end{abstract}
\maketitle

\section{Introduction and statement of results}
In this paper we examine conditions on a $C^*$-dynamical system
$(A,G,\alpha)$ so that induction of primitive ideals from stabilizers
always leads to primitive ideals of the (full) crossed product
$A\rtimes_\alpha G$.  For any $C^*$-algebra $B$, we denote by
$\mathcal I(B)$ the set of all closed two-sided ideals of $B$ equipped
with Fell's topology (see \cite{ech:mams96}*{Chapter~1} for the
definition).  Recall (cf., \cite{rw:morita}*{\S3.3}) that there are
continuous maps
\begin{equation*}
  \Ind_{H}^{G}:\mathcal{I}(A\rtimes_{\alpha}H)\to
  \mathcal{I}(\acg)\quad \text{and}\quad \Res:\mathcal{I}(\acg)\to
  \mathcal{I}(A) 
\end{equation*}
characterized by
\begin{gather}
  \Ind_{H}^{G} \big(\ker (\sigma\rtimes V)\big) =  \ker\big( \Ind_{H}^{G} (\sigma\rtimes V)\big) \\
  \intertext{for any nondegenerate representation $\sigma\rtimes V$ of
    $A\rtimes_{\alpha}H$, and} \Res \big(\ker (\pi\rtimes U) \big)=
  \ker \pi
\end{gather}
for any nondegenerate representation $\pi\rtimes U$ of $\acg$, where
$\Ind$ and $\Res$ denote induction and restriction of representations.
For a primitive ideal $P$ we denote by
\begin{equation*}
  G_{P}:=\set{s\in G:\alpha_{s}(P)=P}
\end{equation*}
the stability group at $P$ for the continuous action of $G$ on
$\Prim(A)$.

\begin{definition}
  \label{def-std-def}
  We say that the dynamical system $\aga$ satisfies the 
    \emph{Effros-Hahn-induction property} (\EHI) if, given $P\in\Prim(A)$
  and a primitive ideal $J$ in $\Prim (\acgp)$ with $\Res J=P$, then
  $\indgpg J$ is a primitive ideal in $\acg$.
  We say $\aga$ satisfies the \emph{strong
    Effros-Hahn-induction property} (strong-\EHI) if, given
  $P\in\Prim(A)$ and an irreducible representation $\rho\rtimes V$ of
  $\acgp$ with $\ker\rho=P$, then $\indgpg (\rho\rtimes V)$ is
  irreducible.
\end{definition}

It is clear that strong-\EHI{} implies \EHI.  It is well-known
that all separable systems $(A,G,\alpha)$ with $A$ \emph{commutative}
satisfy strong-\EHI{}.  The proof goes back to Mackey
(\cite{mac:pnasus49}*{\S6}; see also \cite{gli:pjm62}*{pp.~900-901}).
In fact, the separability assumption can be dropped by
\cite[Proposition 4.2]{wil:tams81}.  Building on
\cite{wil:tams81}*{Proposition~4.2}, Olesen and Raeburn were able to
show that if $A$ is a separable continuous-trace \cs-algebra with $G$
acting freely on the spectrum, then $(A,G,\alpha)$ satisfies
strong-\EHI{} \cite{olerae:jfa90}*{Lemma~3.2}.\footnote{Olesen and
  Raeburn assumed that $G$ was abelian, but their proof was extended  
  to the non-abelian case by Deicke in \cite{deicke}. }  It follows from
Green's proof of \cite{gre:am78}*{Theorem~24} that if $(A,G,\alpha)$
is separable with $G$ amenable and acting freely on $\Prim A$, then
$(A,G,\alpha)$ satisfies \EHI.\footnote{Green's proof requires only
  that $(A,G,\alpha)$ be what he calls ``EH-regular''.  The
  Gootman-Rosenberg-Sauvageot result \cite{gooros:im79} implies this
  for separable systems with $G$ amenable. (Note that $\Ind I = \Ind
  \bigcap_{s}s\cdot I$ by \cite{gre:am78}*{Proposition~11(ii)}.)  The
  result for free actions on separable systems follows from
  Proposition~\vref{prop-normal-stability}.}  Green's result
generalizes a result of Zeller-Meier for $G$ discrete
\cite{zel:jmpa68}*{Theorem~5.15}.  We would like to believe that the
following conjecture is true.
\begin{conjecture}
  \label{conj-1} All separable dynamical systems $\aga$ satisfy
  \EHI.
\end{conjecture}

If we add the hypothesis that $G$ is amenable, then \EHI{} is claimed
as part of the first author's \cite{ech:mams96}*{Theorem~1.4.14}.
However, \cite{ech:mams96}*{Theorem~1.4.14} is meant to be a summary
of the main results in
Gootman-Rosenberg \cite{gooros:im79}, and for the result we are
interested in, Gootman and Rosenberg are relying on
Sauvageot's \cite{sau:jfa79}*{Proposition~2.1}.  The problem
 is that
Sauvageot does not work with an arbitrary primitive ideal in $\acgp$
with restriction $P$.  He requires an additional assumption that
certain associated representations are homogeneous:
\begin{definition}
  \label{def-homogeneous}
  A representation $\pi:A\to B(\H_{\pi})$ is called \emph{homogeneous}
  if every non-zero sub-representation of $\pi$ has the same kernel as
  $\pi$.
\end{definition}

Sauvageot's result is as follows.

\begin{prop}[{\cite[Proposition~2.1]{sau:jfa79}}]
  \label{prop-sauv-2.1}
  Let $(A,G,\alpha)$ be a separable dynamical system.  Suppose that
  $\rho$ is a \emph{homogeneous representation} of $A$ with kernel
  $P$, and that $\rho\rtimes V$ is a \emph{homogeneous representation}
  of $\acgp$.  Then $\indgpg(\rho\rtimes V)$ is a homogeneous
  representation of $\acg$.
\end{prop}

\begin{remark}
  \label{rem-homo-basics} It is clear that every irreducible
  representation is homogeneous.  On the other hand, a representation
  $\pi:A\to B(\H_{\pi})$ is homogeneous if and only if given $I\in
  \mathcal{I}(A)$, then $\overline{\pi(I)\H_{\pi}}$ is either all of
  $\H_{\pi}$ or $\set0$.  This is the content of
  \cite[Lemme~1.5]{sau:jfa79} or \cite[Theorem~1.4]{eff:tams63}.  It
  follows that the kernel of a homogeneous representation is always a
  \emph{prime} ideal \cite[Corollary~1.5]{eff:tams63}.  Since every
  prime ideal is primitive in a \emph{separable} \cs-algebra
  \cite[Theorem~A.49]{rw:morita}, the kernel of any homogeneous
  representation of a separable \cs-algebra is a primitive
  ideal.\footnote{Separability may be crucial here.  Weaver has
    constructed a nonseparable example of a prime ideal which is not
    primitive \cite{wea:jfa03}.}
  \end{remark}
\begin{remark}
   As the emphasis in Proposition \vref{prop-sauv-2.1} indicates,
  our problem is
  resolving Sauvageot's ``extra'' hypothesis on $\rho$ with
  Conjecture~\vref{conj-1}.  Unfortunately, Sauvageot's proof
  makes significant use of the homogeneity of $\rho$ as it
  is necessary to view $\indgpg(\rho\rtimes V)$ as a direct integral
  of homogeneous representations over $G/G_{P}$.  Gootman and
  Rosenberg are careful to note this property in their definition of
  an ``induced primitive ideal'' just prior to stating their
  Theorem~3.1 on page~290 of \cite{gooros:im79}.  Unfortunately, these
  assumptions were neglected by the first author when he included
  \EHI{} as part of
  \cite{ech:mams96}*{Theorem~1.4.14}, so we do not know whether that
  result is valid in the full generality as stated there.
  Fortunately, this does not affect any other results of
  \cite{ech:mams96}, since Propositions \ref{prop-ehi} and
  \ref{prop-normal-stability} stated below clearly imply \EHI{} in all
  situations where \cite{ech:mams96}*{Theorem~1.4.14} was used in that
  \emph{Memoir}.
\end{remark}

Using techniques similar to Sauvageot's, in \S\ref{sec-ind} we prove the
following variation of Proposition~\ref{prop-sauv-2.1}:

\begin{thm}
  \label{thm-2.1-irred}
  Let $(A,G,\alpha)$ be a separable dynamical system.  Suppose that
  $\rho$ is a homogeneous representation of $A$ with kernel $P$, and
  that $\rho\rtimes V$ is an irreducible representation of
  $A\rtimes_{\alpha}G_{P}$.  Then $\Ind_{G_{P}}^{G}(\rho\rtimes V)$ is
  an irreducible representation of $\acg$.
\end{thm}

In order to use the above propositions for proving (strong-)\EHI, it
would certainly be sufficient to show that for any irreducible
representation $\rho\rtimes V$ of $A\rtimes_{\alpha}G_P$ with $\ker
\rho=P$, the representation $\rho$ is automatically homogeneous.
Unfortunately, this turns out to be not true in general (see Example
\ref{ex-conj-false} below).  However, we shall see that it \emph{is}
true if $P\in \Prim(A)$ is \emph{locally closed}; that is, if the
point set $\set P$ is open in $\overline{\set P}$.  There are actually
many $C^*$-algebras with the property that all points in $\Prim(A)$
are locally closed; for example, all type I $C^*$-algebras have this
property.  As a consequence of this and 
Theorem~\ref{thm-2.1-irred} we get

\begin{prop}
  \label{prop-ehi}
  Suppose that $\aga$ is a separable dynamical system such that all
  points are locally closed in $\Prim(A)$.  Then $\aga$ satisfies
  strong-\EHI. In particular, if $\aga$ is a separable dynamical
  system such that $A$ is type~I, then $\aga$ satisfies strong-\EHI.
\end{prop}

Another positive result towards Conjecture~\ref{conj-1} is the
following (see \S\ref{sec-EHI} for the proof):

\begin{prop}
  \label{prop-normal-stability}
  Suppose that $(A,G,\alpha)$ is a separable dynamical system such
  that $G_{P}$ is normal in $G$ for all $P\in\Prim(A)$ (which is
  clearly true if $G$ is abelian).  Then $\aga$ satisfies strong-\EHI.
\end{prop}

If we only consider the weaker property \EHI{} and if we are willing
to assume that all stabilizers are amenable, the normality assumption
on the stabilizers used above can be weakened to obtain the following
result (see \S\ref{sec-EHI}):

\begin{prop}
  \label{prop-normalizer}
  Suppose that $\aga$ is a separable $C^*$-dynamical system such that
  all stabilizers $G_P$ are amenable and such that for all $P,Q\in
  \Prim(A)$ satisfying
$$P=\bigcap_{s\in G_Q}s\cdot Q$$
we have $G_Q\subseteq N(G_P)$ or $G_P\subseteq N(G_Q)$ (where
$N(H)=\{s\in G: sHs^{-1}\subseteq H\}$ denotes the {\em normalizer} of
a subgroup $H$ of $G$).  Then $\aga$ satisfies \EHI.
\end{prop}

A provocative program for producing a counterexample to
Conjecture~\ref{conj-1} or to strong-\EHI, 
would involve finding an example of a group action on a $C^*$-algebra
\emph{which does not} satisfy the normality condition on the
stabilizers as used in the above proposition.  However, this seems to
be difficult.  The importance of settling Conjecture~\ref{conj-1} is
discussed in more detail in \S\ref{sec:problem-with-current}.

 \medskip

In the appendix, the first author uses this opportunity to give a
corrected proof of \cite[Theorem 5.5.13]{ech:mams96}. This is
necessary because the original proof relies on \cite[Lemma
5.5.17]{ech:mams96}, which turned out to be false.

\subsection*{Acknowledgements} 
Much of this work has been done while the first author visited
Dartmouth College in January 2005. He wants to take this opportunity
to thank the members of the Department of Mathematics, and
especially Professor Dana Williams, for their very warm hospitality.

\section{Irreducible induced representations}\label{sec-ind}

In this section we want to give the proof of 
Theorem~\ref{thm-2.1-irred} as stated in the introduction. As mentioned
there, we will make extensive use of the ideas developed by Sauvageot
in \cite{sau:jfa79}.  Let us first review the process of inducing
representations (cf., e.g., \cite{rw:morita}).  Let $H$ be a closed subgroup of $G$ and let
$(\rho,V)$ be a covariant representation of $(A,H,\alpha\restr H)$ on
$\H$.  Let $\X=\overline{C_{c}(G,A)}$ be Green's $\E\sme\ach$-\ib,
where $\lt:G\to \Aut\big(C_0(G/H)\big)$ denotes left translation.  Let
$\xind (\rho\rtimes V)$ be the representation of $E:=\E$ corresponding
to $\rho\rtimes V$ via $\X$. Then $\xind (\rho\rtimes V)$ acts on the
Hilbert space $\V$ which is the completion of $C_{c}(G,A)\odot \H$
with respect to the pre-inner product
\begin{equation*}
  \ip(f\tensor h|g\tensor k):=\bip(\rho\rtimes V(\rip\ach<g,f>)h|k),
\end{equation*}
where
\begin{equation*}
  \rip\ach<g,f>(t)=\gamma(t)\int_{G}\alpha_{s}^{-1}\bigl(g(s)^{*}f(st)
  \bigr) \,d\mu_{G}(s)\quad\text{and}\quad\gamma(t)
  :=\Bigl(\frac{\Delta_{G}(t)}{\Delta_{H}(t)}\Bigr)^{\frac12}.
\end{equation*}
If $c\in E_{0}:=C_{c}(G\times G/H,A)$, then $(\xind (\rho\rtimes
V))(c)$ sends the class of $f\tensor h$ to the class of $c\cdot
f\tensor h$ where
\begin{equation*}
  c\cdot f(s):=\int_{G}c(r,\dot
  s)\alpha_{r}\bigl(f(r^{-1}s)\bigr)\,d\mu_{G}(r) .
\end{equation*}
The representation $\Ind_{H}^{G}(\rho\rtimes V)$ of $\acg$ induced
from $\rho\rtimes V$ is, by definition, the restriction of $\xind
(\rho\rtimes V)$ to the image of $\acg$ in $M(E)$.  Thus, if $g\in
C_c(G,A)\subseteq \acg$, then $\Ind_{H}^{G}(\rho\rtimes V)(g)$ maps
the class of $f\tensor h$ to the class of $g*f\tensor h$.

\subsection{\boldmath Realizing Induced Representations on $\Hind$}
We want to realize both these induced representations on the classical
space $\Hind$ (cf. \cite{fol:course}*{\S6.1}).  As described in
\cite{fol:course}*{\S6.1, Remark~1}, we can
realize $\Hind$ as the collection of $\mu_{G}$-almost everywhere
equivalence classes of functions in $\sHind$ consisting of Borel
functions $\xi:G\to\H$ such that \emph{for all} $s\in G$ and $t\in H$
\begin{equation*}
  \xi(st)=V_{t}^{-1}\bigl(\xi(s)\bigr),
\end{equation*}
and such that
\begin{equation*}
  \int_{G/H}\|\xi(s)\|^2\,d\beta(\dot s)<\infty,
\end{equation*}
where $\beta$ is a quasi-invariant measure on $G/H$ corresponding to a
positive continuous function $\varphi:G\to (0,\infty)$ such that
\begin{equation*}
  \varphi(st)=\frac{\Delta_{H}(t)}{\Delta_{G}(t)}\varphi(s)\quad\text{for
    all $s\in G$ and $t\in H$.}
\end{equation*}

$\Hind$ is a Hilbert space with respect to the inner product
\begin{equation*}
  \ip(\xi|\eta):=\int_{G/H}\bip(\xi(s)|\eta(s))\, d\beta(\dot s).
\end{equation*}

We can define $W:C_{c}(G,A)\odot\H\to\sHind$ by
\begin{equation*}
  W(f\tensor
  h)(r):=\int_{H}\varphi(rt)^{-\frac12}\rho\bigl(\alpha_{r}^{-1}\bigl(
  f(rt)\bigr)\bigr) V_{t}h\,d\mu_{H}(t).
\end{equation*}
Then $W$ extends to a unitary transformation from $\V$ onto $\Hind$.

Any representation of $\E$ is determined by a covariant representation
$(\pi, U)$ on $(A,G,\alpha)$ together with a representation $M$ of
$C_{0}(G/H)$ which commutes with $\pi$.  It is not hard to check that
$W$ intertwines $\xind (\rho\rtimes V)$ with the representation
determined by $\pi$, $U$ and $M$ given by
\begin{align*}
  M(\phi)\xi(s)&:=\phi(\dot s)\xi(s) \\
  \pi(a)\xi(s)&:= \rho\bigl(\alpha_{s}^{-1}(a)\bigr)\xi(s)\\
  U(r)\xi(s)&:=
  \Bigl(\frac{\varphi(r^{-1}s)}{\varphi(s)}\Bigr)^{\frac12}\xi(r^{-1}s).
\end{align*}
The representation intertwined with $\Ind_{H}^{G}(\rho\rtimes V)$ is
given by $\pi\rtimes U$.

\subsection{\boldmath Realizing $\pi$ as a Direct Integral over $G/H$}
We want to see that $\pi$ can be viewed as a very elementary direct
integral.  Actually, we don't need to invoke direct integrals in their
full glory (as described, for example, in \cite{dix:von}).  It will be
enough to rely on the treatment in \cite{arv:invitation}*{\S4.2} which
considers the special case where the underlying Hilbert space is of
the form $L^{2}(X,\mu;\H)$.  To apply these results, we need a few
observations.  Let $c:G/H\to G$ be a Borel cross section for the
quotient map, and let
$b:G\to H$ be the Borel map defined by $b(s)=c(\dot s)^{-1}s$ for all
$s\in G$.  Note that $b(st)=b(s)t$ for all $s\in G$ and $t\in H$.  If
$f\in\L^{2}(G/H,\beta;\H)$, then
\begin{equation*}
  H_{f}(s):=V_{b(s)}^{-1}\bigl(f(\dot s)\bigr)
\end{equation*}
belongs to $\sHind$ and $\|H_{f}\|_{2}=\|f\|_{L^{2}(G/H;\H)}$.
Conversely, if $h\in\sHind$, then
\begin{equation*}
  F_{h}(\dot s):= h\bigl(c(\dot s)\bigr)
\end{equation*}
is in $\L^{2}(G/H,\beta;\H)$.  The maps $f\mapsto H_{f}$ and $h\mapsto
F_{h}$ are inverses and we obtain an isomorphism between $\Hind$ and
$L^{2}(G/H,\beta;\H)$.  Thus we can view $\sHind$ as the
$L^{2}$-sections of the trivial Borel Hilbert bundle $G/H\times\H$.
By definition, the diagonal operators $\Delta(G/H\times\H,\beta)$ are
naturally identified with $L^{\infty}(G/H)$; if $\psi$ is a bounded
Borel function on $G/H$ then the corresponding diagonal operator on
$\Hind$ is given by extending $M$ in the obvious way:
\begin{equation*}
  M(\psi)\xi(s):=\psi(\dot s)\xi(s).
\end{equation*}
As shown in \cite{arv:invitation}*{Theorem~4.2.1}, the decomposable
operators are exactly the commutant $M(L^{\infty}(G/H))'$, and they are
identified with equivalence classes of bounded weak-operator Borel
functions $F:G/H\to B(\H)$ via
\begin{equation*}
  T_{F}(h)(s):= 
  \widetilde F(s)\bigl(h(s)\bigr), \quad\text{where}\quad\widetilde
  F(s):= V_{b(s)}^{-1}F(\dot s)V_{b(s)}\quad\text{and}\quad h\in
  \L_V^2(G,\beta;\H). 
\end{equation*} 
Note that $\tilde F$ satisfies
\begin{equation}
  \label{eq:333}
  \widetilde F(st)=V_{t}^{-1}\widetilde F(s)V_{t}\quad\text{for all
    $s\in G$ and $t\in H$.}
\end{equation}
If $\widetilde F:G\to B(\H)$ is any weak-operator Borel function
satisfying \eqref{eq:333}, then $F(\dot s):=\widetilde F\bigl(c(\dot
s)\bigr)$ is a weak-operator Borel function on $G/H$ such that
\begin{equation*}
  T_{F}h(s)=\widetilde F(s)\bigl(h(s)\bigr).
\end{equation*}
Thus we can identify the decomposable operators on $\Hind$ with the
bounded weak-operator Borel functions $F:G\to B(\H)$ that transform as
in \eqref{eq:333}. Since
\begin{equation*}
  \rho\bigl(\alpha_{st}^{-1}(a)\bigr)=
  V_{t}^{-1}\rho\bigl(\alpha_{s}^{-1}(a)\bigr) V_{t},
\end{equation*}
$s\mapsto \rho\bigl(\alpha_{s}^{-1}(a)\bigr)$ is a bounded
weak-operator Borel function transforming as in \eqref{eq:333}.
Therefore we can express $\pi$ as a direct integral
\begin{equation}\label{eq:14}
  \pi =\int_{G/H}^{\oplus} \rho_{\dot s}\,d\beta^{P}(\dot s),
\end{equation}
where
\begin{equation}\label{eq:377}
  \rho_{\dot s}(a)h(\dot s)= \rho\bigl(\alpha_{c(\dot
    s)}^{-1}(a)\bigr)\bigl(h(\dot s)\bigr).  
\end{equation}

\subsection{\boldmath The Role Played by the Homogeneity of $\pi$}
From here on, we shall ignore the unitary $W$ and identify
$\Ind_{H}^{G}(\rho\rtimes V)$ with $(\pi, U)$ and $\xind (\rho\rtimes V)$ with
$\bigl((M\tensor\pi), U\bigr) $.  We will also specialize to the
situation where $H=G_{P}$ and $(\rho,V)$ is an \emph{irreducible}
representation with $\rho$ a \emph{homogeneous
representation} with kernel $P$.  This assumption on $\rho$ is used to
show that all operators in the commutant of $\pi$ are decomposable.
This is the essential step in both our proof of
Theorem~\vref{thm-2.1-irred} and Sauvageot's proof of
Proposition~\vref{prop-sauv-2.1}.

\begin{proof}[Proof of Theorem~\vref{thm-2.1-irred}]
Decompose $\pi$ as in \eqref{eq:14}.  Since $\rho$ is homogeneous with
kernel $P$, each $\rho_{\dot s}$ is homogeneous with kernel $s\cdot P$.
Thus if $\dot s\not=\dot r$, then $\ker\rho_{\dot s}\not=\ker\rho_{\dot
  r}$. Thus \eqref{eq:377} is an example of the sort of decomposition
studied by Effros in \cite{eff:tams63}, and what is now called an
ideal center decomposition. One of the major theorems of Effros's
theory is that for such a decomposition, the diagonal operators
$$M\bigl(L^{\infty}(G/G_P)\bigr)=\Delta(G/G_{P}\times\H,\beta)$$ coincide with
a certain von-Neumann subalgebra, $\IC(\pi)$, of the center
$\pi(A)'' \cap \pi(A)'$
(\cite{eff:tams63}*{Theorem~1.10}).  Since
$M\bigl(L^{\infty}(G/G_P)\bigr)=\Delta(G/G_{P}\times\H,\beta)=\IC(\pi)\subset
\pi(A)'' \cap \pi(A)'$, we have
\begin{equation}
  \label{eq:13}
  \pi(A)'\subset M\bigl(L^{\infty}(G/G_{P})\bigr)' 
\subset M\bigl(C_{0}(G/G_{P})\bigr)'. 
\end{equation}
Thus if $T$ is in the commutant of $\Ind_{G_{P}}^{G}(\rho\rtimes V)=\pi\rtimes U$, then $T\in
U(G)'$ and
$T\in\pi(A)'$. It follows then from \eqref{eq:13} that $T\in
M\bigl(C_{0}(G/G_{P})\bigr)'$ too.  But then $T$ is in the
commutant of the irreducible representation $\xind (\rho\rtimes V)=(M\otimes\pi)\rtimes U$.  
Hence $T$ is a
scalar operator and the commutant of $\Ind_{G_{P}}^{G}(\rho\rtimes V)$ is trivial. 
Thus $\Ind_{G_{P}}^{G}(\rho\rtimes V)$ is irreducible. 
\end{proof}

\section{(Strong-)\EHI{} and examples}\label{sec-EHI}

In this section we want to deduce all remaining results as stated in the introduction
and present some examples.
We start with recording the following elementary consequence
of the definition of the Jacobson Topology on $\Prim(A)$.
 
\begin{lemma}\label{lem-popen}
  Suppose that $A$ is a \cs-algebra and that $P\in\Prim(A)$.  Then
  $\set P$ is open in $\overline{\set P}$ if and only if
$
    J:= \bigcap\{P'\in\Prim (A):  P'\supset P ,\, P'\not= P\}
$
  \emph{properly} contains $P$. (By convention, the intersection over the empty set is the whole
  algebra $A$)
\end{lemma}
\begin{proof}
  If $\set P$ is open in its closure, then there is an ideal
  $I\in\mathcal{I}(A)$ such that
  \begin{equation*}
    \mathcal{O}_{I}:= \set{P'\in\Prim(A) : P'\not\supset I}
  \end{equation*}
  contains $P$ and is disjoint from $\overline{\set P}\setminus \set
  P$.  But then it follows that $J\supset I$, and $P\not \supset I$
  forces $J\not=P$ as required.
  Conversely, if $J$ properly contains $P$, then $\mathcal{O}_{J}=\{P'\in \Prim(A): P'\not\supset J\}$ is
  an open set containing $P$ and disjoint from $\overline{\set
    P}\setminus \set P$.
\end{proof}

We use this result to show

\begin{lemma}
\label{lem-locally-closed}
Let $\aga$ be a $C^*$-dynamical system and assume that $P\in\Prim(A)$ is locally closed.
Suppose that $\rho\rtimes V$ is an irreducible representation of $\acgp$
  with $\ker \rho = P$. Then $\rho$ is homogeneous.
 \end{lemma}
\begin{proof}
   In view of
  Remark~\vref{rem-homo-basics}, it will suffice to see that if $I$ is
  an ideal in $A$ such that $I\not\subset P$, then
  \begin{equation*}
    \overline{\rho(I)\H_{\rho}}=\H_{\rho}. 
  \end{equation*}
  Since $\rho(I)=\rho(I+P)$, we can assume that $I$ strictly contains
  $P$.  If $s\in G_{P}$, then $\alpha_{s}(I)$ strictly contains $P$,
  and since every ideal in $A$ is the intersection of the primitive
  ideals which contain it, our assumption on $\set P$ and
  Lemma~\ref{lem-popen} imply that
\begin{equation*}
  K:= \bigcap_{s\in G_{P}}\alpha_{s}(I) \supset
  \bigcap_{\substack{P'\in\Prim(A)\\ P'\supset P \\ P'\not= P}} P':= J,
\end{equation*}
properly contains $P$.  But $K$ is a $G_{P}$-invariant ideal so that
$
  \overline{\rho(K)\H_{\rho}}
$
is nonzero and invariant for both $\rho$ and $V$: if $a\in K$ and
$s\in G_{P}$, then
\begin{equation*}
  V_{s}\rho(a)h = \rho\bigl(\alpha_{s}(a)\bigr)V_{s}h \in
  \rho(K)\H_{\rho}. 
\end{equation*}
Since $\rho\rtimes V$ is irreducible by assumption,
$\overline{\rho(K)\H_{\rho}}$ must be all of $\H_{\rho}$.  Since
$\rho(I)\H_{\rho} \supset \rho(K)\H_{\rho}$, this suffices. 
\end{proof}

\begin{proof}
  [Proof of Proposition~\vref{prop-ehi}] The result will follow
  directly from Theorem~\vref{thm-2.1-irred} together with Lemma
  \ref{lem-locally-closed} as soon as we see that primitive ideals of
  type I algebras are locally closed. We can identify $\overline{\set
    P}$ with $\Prim(A/P)$ so that we may as well assume that $A$ is
  primitive (i.e., $\set 0$ is a primitive ideal) and prove that $\set
  {\{0\}}$ is open in $\Prim(A)$.  Since $A$ is type~I, it follows
  that $A$ must contain an essential ideal $I$ isomorphic to the
  compact operators on some Hilbert space.  In particular, every
  nonzero ideal in $A$ contains $I$.  Now the result follows from
  Lemma~\ref{lem-popen}.
\end{proof}

\begin{example}
  \label{ex-bad}
There are separable \cs-algebras with primitive ideals
  $P$ such that $\set P$ is not open in $\overline{\set P}$.  To see an example, let $A$
  be the \cs-algebra $\cs(E)$ associated to the graph $E$ with
  vertices $E^{0}=\Z^{2}$, and with edges connecting $(i,j)$ to both
  $(i+1,j)$ and $(i,j+1)$:
  
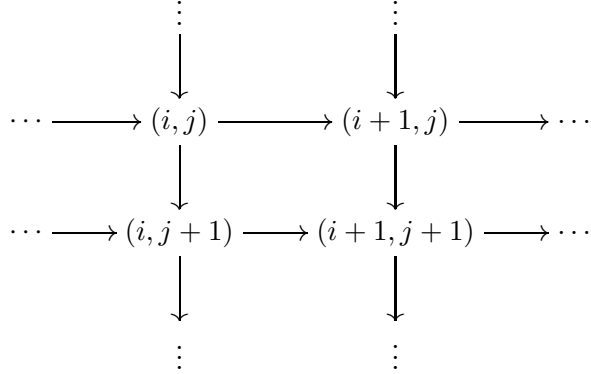
\begin{equation*}
    \xymatrix{&\protect\vdots\ar[d]&\protect\protect\vdots\ar[d]\\
   \cdots\ar[r]&(i,j)\ar[r]\ar[d]&
    (i+1,j)\ar[r]\ar[d]&\dots \\
    \cdots \ar[r]&(i,j+1)\ar[r] \ar[d]&
    (i+1,j+1) \ar[r]\ar[d]&\dots \\
    &\protect\vdots&\protect\vdots} 
  \end{equation*}
 We refer to \cite{bprs:nyjm00} for the details of the construction of $C^*(E)$.
  (It was Astrid an Huef who pointed out to us that this algebra had an
  ``interesting'' primitive ideal structure.) 

Since $E$ has no loops, it satisfies condition~(K)  of \cite{kprr:jfa97}, and
\cite{bprs:nyjm00}*{Proposition~6.1} implies that the primitive ideals of
$\cs(E)$ are in one-to-one correspondence with the set of maximal
tails in $E$.  Recall that if $\gamma\subset E^{0}$ is a maximal tail,
then the corresponding primitive ideal is $I_{H}$ which is associated
to the saturated hereditary subset $H:=E^{0}\setminus \gamma$.  If $H$
is any saturated hereditary subset in $E^{0}$ (not necessarily
corresponding to a maximal tail), then $I_{H}$ is the ideal generated
by the projections $p_{v}$ with $v\in H$.  Furthermore, $p_{v}\in
I_{H}$ if and only if $v\in H$.  (For references and further
discussion, see \cite[\S4]{bprs:nyjm00}.)  It is not hard to see that the
maximal tails in $E$ are $E^{0}$ itself together with
\begin{equation*}
  \gamma_{n}:=\set{(i,j):j\ge n}\quad\text{and}\quad \zeta_{m}:=
  \set{(i,j): i\ge m}. 
\end{equation*}
In particular, $\gamma_{n}$ corresponds to the primitive ideal
$I_{H_{n}}$ associated to
\begin{equation*}
  H_{n}:=\set{(i,j):j<n},
\end{equation*}
and $E^{0}$ corresponds to the zero ideal.  (So $\cs(E)$ is
primitive.)  Now
\begin{equation*}
  I:= \bigcap_{n=1}^{\infty}I_{H_{-n}}
\end{equation*}
is an ideal which contains none of the projections $p_{(i,j)} $. 
Thus, $I=\set0$ and Lemma~\vref{lem-popen} implies that $\set {\{0\}}$
is not open in $\Prim(A)$. 
\end{example}

Indeed, the above example may be extended to give an example 
of a $C^*$-dynamical system $\aga$, a primitive ideal $P\in \Prim(A)$, and 
an irreducible representation $\rho\rtimes V$ of $\acgp$ with $\ker\rho=P$,
such that $\rho$ \emph{is  not homogeneous}:
\begin{example}
\label{ex-conj-false}
Let $A=\cs(E)$ where $E$ is the directed graph of
Example~\ref{ex-bad}.  The universal property of graph algebras
guarantees that there is an automorphism $\phi$ of $\cs(E)$ which
takes the projection $p_{(i,j)}$ to $p_{(i+1,j)}$.  Thus in the
notation of Example~\ref{ex-bad}, we have
\begin{equation*}
  \phi(I_{H_{n}})=I_{H_{n+1}}.
\end{equation*}
Let
\begin{equation*}
  Q:=I_{H_{0}}\quad\text{and}\quad P:=\set 0=\bigcap_{n\in\Z}n\cdot Q.
\end{equation*}
Then $P$ and $Q$ are primitive ideals in $\cs(E)$ and we have
\begin{equation*}
  G_{Q}=\set e\quad\text{and}\quad G_{P}=\Z.
\end{equation*}
Let $\rho\in\specnp{\cs(E)}$ be such that $\ker\rho=Q$.  Then
Theorem~\vref{thm-2.1-irred} implies that
\begin{equation*}
  \Ind_{\set e}^{G}\rho
\end{equation*}
is an irreducible representation of
$\cs(E)\rtimes_{\phi}\Z$.\footnote{Note that by
  \cite{kumpas:etds99}*{Theorem~1.1}, $\cs(E)\rtimes_{\phi}\Z\cong
  \K$.}  Since $G_{P}=\Z$, it follows that $
  \pi\rtimes U:=\Ind_{\set e}^{G}\rho$
is an irreducible representation of $\cs(E)\rtimes_{\phi}G_{P}$ such that   $\pi\cong\bigoplus_{n\in\Z}
(\phi^n\circ \rho)$ and 
$$\ker\pi=\bigcap_{n\in \Z} \ker(\phi^n\circ \rho)=\bigcap_{n\in \Z}n\cdot Q=\set 0=P.$$
Since $\rho$ is a subrepresentation of $\pi$ with $\ker\rho=Q\neq P=\ker\pi$
it follows that $\pi$ is not homogeneous.  
\end{example}

One might expect that the failure of homogeneity of $\rho$ for an
irreducible representation $\rho\rtimes V$ of $\acgp$ with
$\ker\rho=P$ might also lead to a counterexample for (strong-)\EHI.
This idea certainly doesn't work in the above example, since there we
have $G_P=G$. Indeed, Proposition \vref{prop-normal-stability}, which
we are now going to prove, shows that in order to get a counterexample
to strong-\EHI, one has to consider actions of non-abelian
groups:

\begin{proof}[Proof of Proposition \vref{prop-normal-stability}]
  Suppose that $P\in\Prim(A)$ and that $\rho\rtimes V$ is an
  irreducible representation of $A\rtimes_{\alpha}G_{P}$ with $\ker
  \rho=P$.  Since $G_{P}$ is normal in $G$, it follows from
  \cite{gre:am78}*{Proposition 1} (see also
  \cite{ech:mams96}*{Example~1.1.1}) that there is a natural
  isomorphism
  \begin{equation}\label{eq:18}
    A\rtimes_{\alpha}G\cong
    \bigl(A\rtimes_{\alpha}G_{P}\bigr)\rtimes^{\tau}_{\gamma}G 
  \end{equation}
for a twisted system $(A\rtimes_{\alpha}G_{P},G,\gamma,\tau)$ where
(following the notation in \cite{ech:mams96})
\begin{equation*}
  \gamma_{s}(f)(t):=
 \delta(s)\alpha_{s}\bigl(f(s^{-1}ts)\bigr)
\end{equation*}
for $f\in C_{c}(G_{P},A)$,
\begin{equation*}
  \int_{G_{P}}
  g(sts^{-1})\,d\mu_{G_{P}}(t)=\delta(s)\int_{G_{P}}g(t)\,d\mu_{G_{P}}
  (t)
\end{equation*}
for $g\in C_{c}(G_{P})$, and
\begin{equation*}
  \tau(t)(f)(s)=\alpha_{t}\bigl(f(t^{-1}s)\bigr)
\end{equation*}
for $f\in C_{c}(G_{P},A)$ and $t\in G_{P}$.  Similarly, we have an isomorphism
\begin{equation*}
  A\rtimes_{\alpha}G_{P}\cong
  \bigl(A\rtimes_{\alpha}G_{P}\bigr)\rtimes_{\gamma}^{\tau}G_{P}, 
\end{equation*}
and this isomorphism intertwines $\rho\rtimes V$ with $(\rho\rtimes
V)\rtimes V$.  In particular, $(\rho\rtimes V)\rtimes V$ is
irreducible.  Furthermore, \cite{ech:mams96}*{Proposition~1.4.3}
implies that $\Ind_{G_{P}}^{G}\bigl((\rho\rtimes V)\rtimes V\bigr)$ is
intertwined with $\Ind_{G_{P}}^{G}(\rho\rtimes V)$ by the isomorphism
of \eqref{eq:18}.  Therefore it will suffice to show that
\begin{equation}
  \label{eq:19}
  \Ind_{G_{P}}^{G}\bigl((\rho\rtimes V)\rtimes V\bigr)
\end{equation}
is irreducible.  To prove this, we want to apply
Theorem~\vref{thm-2.1-irred} to the dynamical system
$(A\rtimes_{\alpha}G_{P},G,\gamma)$ and the irreducible representation
$(\rho\rtimes V)\rtimes V$.
Since $\rho\rtimes V$ is
irreducible, and therefore homogeneous, we only need to check that
$G_{P}$ is equal to $G_{J}$ where $J$ is equal to $\ker(\rho\times V)$
and $G_{J}$ is the stability group of
$J\in\Prim(A\rtimes_{\alpha}G_{P})$ with respect to the action induced
by $\gamma$.  However, as $\gamma$ admits a twist $\tau$ with
$N_{\tau}=G_{P}$, we have $G_{P}\subset G_{J}$.  On the other hand, it
is not hard to check that
\begin{equation*}
  (\rho\rtimes V)\circ \gamma_{s}^{-1} = s\cdot \rho\rtimes {}^{s}V,
\end{equation*}
where, as usual, $s\cdot \rho:=\rho\circ\alpha_{s}^{-1}$ and
${}^{s}V(t):=V(sts^{-1})$.  In particular,
\begin{equation*}
  \Res\bigl(\gamma_{s}(J)\bigr)=\alpha_{s}(P).
\end{equation*}
Thus $\gamma_{s}(J)\not=J$ if $s\notin G_{P}$.  This completes the proof.
\end{proof}

\subsection{A Method in Search of a Counterexample}

Here we see that using the positive solution of the Effros-Hahn
conjecture for amenable groups, the ``Gootman-Rosenberg-Sauvageot-Theorem''
(\cite[Theorem 3.1]{gooros:im79}), we can obtain a
quasi-result. We are not sure whether this result should be viewed as
evidence to support Conjecture~\vref{conj-1}, or whether it indicates
we should search for
a counter\-example.  

Anyway, suppose that $\aga$ is a separable dynamical system and that
for each $P\in\Prim(A)$, $G_{P}$ is amenable.  Fix
$J\in\Prim\bigl(A\rtimes_{\alpha}G_{P}\bigr)$ with $\Res J=P$.  Since
$G_{P}$ is amenable, we can apply the GRS-Theorem to
$(A,G_{P},\alpha)$, and conclude that there is a $Q\in\Prim(A)$ and a
$K\in \Prim\bigl(A\rtimes_{\alpha}(G_{P})_{Q}\bigr)$ such that
$K=\ker(\sigma \rtimes W)$ with both $\sigma\rtimes W$ and $\sigma$ homogeneous,
and with $\ker\sigma=Q$, and
such that $\Ind_{(G_{P})_{Q}}^{G_{P}}K=J$.  Notice that we must have
\begin{equation}\label{eq:22}
  P=\bigcap_{s\in G_{P}}s\cdot Q.
\end{equation}
Furthermore,
\begin{equation*}
  (G_{P})_{Q}=G_{Q}\cap G_{P}.
\end{equation*}
Now \emph{if} $G_{Q}\subset G_{P}$, then $(G_{P})_{Q}=G_{Q}$
and we have
\begin{align}
  \Ind_{G_{P}}^{G}J&= \Ind_{G_{P}}^{G}
  \bigl(\ker\bigl(\Ind_{G_{Q}}^{G_{P}}\sigma\rtimes W\bigr)\bigr) \notag\\
&= \ker\bigl(\Ind_{G_{Q}}^{G}\sigma\rtimes W\bigr)\label{eq:21}
\end{align}
It follows that \eqref{eq:21} is primitive by Sauvageot's
Proposition~\vref{prop-sauv-2.1}.  Thus we obtain the following
awkwardly stated result:
\begin{lemma}
  \label{lem-awk}
  Suppose that $\aga$ is a separable dynamical system such that
  $G_{P}$ is amenable for all $P\in\Prim(A)$.  Suppose that for all
  $P,Q\in \Prim(A)$ such that \eqref{eq:22} holds, we have
  $G_{Q}\subset G_{P}$.  Then $\aga$ satisfies \EHI.
\end{lemma}

\begin{remark}
  \label{rem-normalizers}
  If \eqref{eq:22} holds, then we have $G_{Q}\subset G_{P}$ provided
  \emph{either}
  \begin{equation*}
    G_{Q}\subset N(G_{P})\quad\text{or}\quad G_{P}\subset N(G_{Q}),
  \end{equation*}
where
\begin{equation*}
  N(H)=\set{s\in G:sHs^{-1}=H}.
\end{equation*}
So the above lemma clearly implies Proposition \vref{prop-normalizer}.
\end{remark}

\section{The Problem with the Current State of the Art}
\label{sec:problem-with-current}

The positive solution of the Effros-Hahn conjecture due to
Gootman-Rosenberg \cite{gooros:im79},  building on the work of Sauvageot
\citelist{\cite{sau:ma77}\cite{sau:jfa79}}, is a critical ingredient in
describing the fine ideal structure of a separable crossed product by
an amenable group.  Nevertheless, the current state of the art leaves
something to be desired.  In order to describe the fine ideal
structure of a \cs-algebra $A$, there are two basic steps.  First it is
necessary to describe $\hat A$ or $\Prim(A)$ as a set.  Secondly, we
want to describe the hull-kernel topology.  Both steps can be formidable.
For many \cs-algebras, especially non type~I \cs-algebras, describing
$\hat A$ can be impossible, and $\Prim(A)$ nearly so.  The general
procedure for
describing either space is to first describe a relatively nice set $X$
consisting of representations --- often called \emph{concrete}
representations --- and then exhibit either $\hat A$ or $\Prim(A)$ as a
quotient of $X$.  One approach to the second step is to equip $X$ with
a topology and show that the hull-kernel topology is the quotient
topology.   

In the case of crossed products $\acg$, a useful example is the case
where $G$ is abelian and $A=C_{0}(\hA)$.  Then a variation on
Theorem~\vref{thm-2.1-irred} produces a continuous map on $\hA\times\widehat
G$ into $\Prim(\acg)$.  The hard bit is to see that this map is
surjective, but his is implied  by the
Gootman-Rosenberg-Sauvageot result.  It is shown in \cite{wil:tams81}
that this map factors through the quotient
\begin{equation}
  \label{eq:15}
  ( \hA\times \widehat G)/\!\!\sim,
\end{equation}
where $(x,\sigma)\sim(y,\rho)$ if
$\overline{G\cdot x}=\overline{G\cdot y}$ and $\sigma\bar\rho\in
G_{x}^{\perp}=G_{y}^{\perp}$, and that we can identify $\Prim
\bigl(C_{0}(\hA)\rtimes G\bigr)$ with \eqref{eq:15} as a
topological space \cite{wil:tams81}*{Theorem~5.3}.

If $A$ is CCR, then we can at least identify a candidate for 
the set $X$ of concrete representations  as follows.
Given $P\in\Prim(A)$, $A/P$ is simple and we can take \emph{any}
irreducible representation of $A/P\rtimes_{\alpha}G_{P}$, lift it
to an irreducible representation $L_{P}$ of $A\rtimes_{\alpha}G_{P}$
with the property that $\Res (\ker L_{P})=P$.  Thus by
Proposition~\ref{prop-ehi}, $\Ind_{G_{P}}^{G}(\ker L_{P})$ is
primitive.  If $G$ is amenable, then this gives us a
map of
\begin{equation}
  \label{eq:16}
  \bigcup_{P\in\Prim(A)}\Prim\bigl(A/P\rtimes_{\alpha}G_{P}\bigr)
\end{equation}
into $\Prim(\acg)$ which is surjective by the GRS-Theorem.

If $A$ is merely GCR, Proposition~\ref{prop-ehi} still applies, but
the picture gets just a bit more complicated.  For any given
$P\in\Prim(A)$, $A/P$ contains a simple essential ideal $K(P)$ which
is an elementary \cs-algebra.  \emph{Any} irreducible representation
$L$ of $K(P)\rtimes_{\alpha}G_{P}$ has a canonical extension $\bar
L=\rho\rtimes V$ to $A/P\rtimes_{\alpha}G_{P}$ with $\rho$ faithful.
Thus we can replace \eqref{eq:16} by
\begin{equation}
  \label{eq:17}
\bigcup_{P\in\Prim(A)}\Prim\bigl(K(P)\rtimes_{\alpha}G_{P}\bigr).
\end{equation}
So even in the GCR case, we have an ``$X$'' from which to start.  One
might even be able to endow $X$ with some sort of topology via a Fell
subgroup algebra construction or some groupoid variation.  A much more
sophisticated and cleaner description is given in
\cite{ech:mams96}*{Theorem~3.1.7 and Remark~3.1.8}. 

In the general case, the additional homogeneity hypotheses
puts some serious holes
in our approach.  Even if $G$ is amenable, given $P\in\Prim(A)$, we
do not see any reason why there has to be a $J\in\Prim
(A\rtimes_{\alpha}G_{P})$ with $\Res J=P$.  All the existing theory
guarantees is that given $K\in\Prim(\acg)$, then there is a $P\in\Prim(A)$ 
and a $J\in\Prim (A\rtimes_{\alpha}G_{P})$ such that $\Res J=P$,
\emph{and of course}, $\Ind_{G_{P}}^{G}J=K$.  Knowing this is a
remarkably powerful tool and certainly gives us plenty of information
about $\Prim(\acg)$.  It is, for example, very useful in determining
the simplicity of $\acg$.  But it definitely does not provide us with
a nice set $X$ in the spirit of \eqref{eq:15} or \eqref{eq:16}.

\section{Appendix: A corrected proof of \cite[Theorem
  5.5.13]{ech:mams96}}

In this appendix, the first author wants to fill a gap in the proof of
Theorem 5.5.13 of the \emph{Memoir} \cite{ech:mams96}, one of the main
results of that \emph{Memoir}. The gap is due to a false lemma
(\cite[Lemma 5.5.17]{ech:mams96}) which was used in the original
proof.  Below, as in \cite{ech:mams96}, $\frak K(G)$ denotes the set
of closed subgroups of $G$ equipped with Fell's topology. $\frak K(G)$
is a compact Hausdorff space and a base of the topology is given by
the sets
$$U(\mathcal F, C)=\{H\in \frak K(G); H\cap V\neq\emptyset\;
\text{for all}\; V\in \mathcal F\;\text{and}\; H\cap C=\emptyset\},$$
where $\mathcal F$ runs through all finite families of open subsets of
$G$ and $C$ runs through the compact subsets of $G$.

In \cite[Lemma 5.5.17]{ech:mams96} we stated that if $N$ is an open
normal subgroup of $G$ and if $q:G\to G/N$ denotes the quotient map,
then the map $H\to q(H); \frak K(G)\to \frak K(G/N)$ is continuous.
Unfortunately, {\em this is not true} in general. To see a
counterexample let $G=\Z$ and $N=2\Z$. Then, if $H_n=(2n+1)\Z$ we have
$H_n\to \{0\}$ in $\frak K(\Z)$ but $q(H_n)=\Z/2\Z\to \Z/2\Z$ in
$\frak K(\Z/2\Z)$. A similar example shows that (contrary to the
statement of the lemma), if $K$ is a compact subgroup of an abelian
group $G$, then the intersection map $H\to H\cap K;\frak K(G)\to \frak
K(K)$ is {\em not continuous} in general.

In what follows we want to correct the proof of \cite[Theorem
5.5.13]{ech:mams96}.  Note that this was the only place in the \emph{Memoir}
where the lemma was used!  \medskip

Recall that if $G$ is an abelian locally compact group and $[\om]\in
H^2(G,\TT)$, then $[\om]$ determines a homomorphism $h_{\om}:G\to
\widehat{G}$ by $h_{\om}(s)(t)=\om(s,t)\om(t,s)^{-1}$.  The group
kernel $\Sigma_{\om}$ of $h_{\om}$ is called the {\em symmetrizer} of
$[\om]$. By a result of Baggett and Kleppner it is known that the
twisted group algebra $C^*(G,\om)$ is type I if and only if $h_{\om}$
has closed range and is open as a map onto its image, in which case it
factors through an isomorphism between $G/\Sigma_{\om}$ and
$\Sigma_{\om}^{\perp}$. We then say that $[\om]$ is {\em type I} (see
the \emph{Memoir} for more details).  As indicated in the original proof of
\cite[Theorem 5.5.13]{ech:mams96}, by passing to a Morita equivalent
untwisted action and by the definition of the Mackey obstruction map
$h^{\alpha,\tau}$ as given in \cite[Definition 5.5.1]{ech:mams96}, the
proof of \cite[Theorem 5.5.13]{ech:mams96} will follow from

\begin{thm}\label{thm}
  Suppose that $(A,G,\alpha)$ is a separable $C^*$-dynamical system
  such that $A$ has continuous trace and $G$ is a compactly generated
  abelian Lie group.  Suppose further that $G$ acts trivially on
  $\hat{A}$ and for each $x\in \Om:=\hat{A}$ let $[\om_x]\in
  H^2(G,\TT)$ denote the Mackey obstruction at $x$ for extending the
  corresponding representation $\rho_x$ to a covariant representation
  of $(A,G,\alpha)$. Let
$$h^{\alpha}:\Om\times G\to \Om \times \widehat{G},\;
h^{\alpha}(x,s)=(x,h_{\om_x}(s)).$$ Then $h^{\alpha}$ has closed range
and is open as a map onto its image if and only if all Mackey
obstructions $[\om_x]$ are type I and the symmetrizer map $\Om\to
\frak K(G)$, $x\mapsto \Sigma_x:=\Sigma_{\om_x}$ is continuous.
\end{thm}

As a consequence of this result (i.e. \cite[Theorem
5.5.13]{ech:mams96}) and \cite[Theorem 5.5.2]{ech:mams96} it follows
that if $(A,G,\alpha)$ is as above, then the crossed product
$A\rtimes_{\alpha}G$ has continuous trace if and only if all Mackey
obstructions are type I and the symmetrizer map is continuous (but see
\cite[Corollary 5.5.14]{ech:mams96} for a more general consequence).

We will need the following basic result which was also partly stated
in \cite[Lemma 5.5.17]{ech:mams96}.

\begin{lemma}\label{lemcorrect}
  Suppose that $G$ is a locally compact group.  If $N$ is an open
  subgroup of $G$ then $\cap_N:\frak K(G)\to \frak K(N);H\mapsto H\cap
  N$ is continuous. Moreover, if $N$ is any closed normal subgroup of
  $G$, and if $\frak K(G)_N$ denotes the set $\{ H\in \frak K(G):
  N\subseteq H\}$, then $q_*:\frak K(G)_N\to \frak K(G/N); H\mapsto
  q(H)$ is a homeomorphism.
\end{lemma}
\begin{proof} Let $H\in \frak K(G)$, let $U_1,\ldots, U_l$ be open
  subsets of $N$, and let $C\subseteq N$ be compact. Since $N$ is open
  in $G$, the sets $U_1,\ldots, U_l$ are also open in $G$.  If we
  denote by $W$ the neighborhood of $H$ in $\frak K(G)$ and by $V$ the
  neighborhood of $H\cap N$ in $\frak K(N)$ determined by these sets,
  then $\cap_N(W)\subseteq V$, which proves continuity of $\cap_N$.

  To see the second assertion let $U_1,\ldots U_l$ be open subsets of
  $G/N$, and let $C\subseteq G/N$ be compact. Choose a compact subset
  $K$ of $G$ such that $q(K)=C$.  Then
$$q_*\left(U\big(q^{-1}(U_1),\ldots, 
  q^{-1}(U_l),K\big)\cap\frak K(G)_N\right)=U\big(U_1,\ldots,
U_l,C\big).$$ Since all elements of $\frak K(G)_N$ contain $N$ as a
subgroup, every basic open set of $\frak K(G)_N$ can be written as
above. Thus, since $q_*$ is bijective, the result follows.
\end{proof}

Recall that every closed subgroup $H$ of a direct product $V\times Z$
of a vector group $V$ with a finitely generated free abelian group $Z$
can be written as a direct product $H= W \times S$ of a vector
subgroup $W$ of $V$ and a subgroup $S$ of $V\times Z$ which is a
finitely generated free abelian group. Thus, for any closed subgroup
$H$ of $V\times Z$ there are two characteristic positive integers,
namely the dimension, $\dim H:=\dim W$, and the rank, $\rank H:=\rank
S$.  If $V_0$ denotes the linear hull of $H$ in $V$, then $\dim V_0=
\dim H+\rank H$. The following is \cite[Lemma 5.5.18]{ech:mams96}.

\begin{lemma}\label{lemrank}
  Let $V$ be a vector group, and let $H_n\to H$ in $\frak K(V)$.  If
  $\rank H\geq \rank H_n$ for all $n\in\NN$, then there exists a
  subsequence $(H_{n_k})_{k\in\NN}$ of $(H_n)_{n\in\NN}$ and elements
  $c_k\in \GL(V)$ such that $c_k\to \id_V$ and $c_k(H)=H_{n_k}$ for
  all $k\in\NN$.
\end{lemma}

If $H$ is a compactly generated abelian Lie group, then $H$ has a (non
canonical) splitting $H= W\times T\times S\times F$, where $W$ is a
vector group (i.e., isomorphic to some $\RR^n$), $T$ is a finite
dimensional torus group, $S$ is finitely generated free abelian, and
$F$ is finite.  Although the splitting is non canonical, the values
$\dim V, \dim T, \rank Z$, and the (isomorphism class of the) finite
group $F$ form a complete invariant of the isomorphism class of $H$.
We need

\begin{lemma}\label{lemsplitt}
  Let $\Sigma$ be a closed subgroup of the direct product $V\times Z$
  of some vector group $V$ with some finitely generated free abelian
  group $Z$. Let $q:V\times Z\to Z$ denote the quotient map and let
  $W\times T\times S\times F$ be a splitting of the compactly
  generated Lie group $G/\Sigma$ as above. Then $\dim T=\rank(V\cap
  \Sigma)$ and $\rank S= \rank Z-\rank q(\Sigma)$. In particular, if
  $G/\Sigma$ is isomorphic to $\Sigma^{\perp}\cong (G/\Sigma)\dach$,
  then $\rank(V\cap\Sigma)=\rank Z-\rank q(\Sigma)$.
\end{lemma}
\begin{proof} Since $V\cdot\Sigma$ is open in $G$ and $q:V\to
  V\cdot\Sigma/\Sigma$ factors through an isomorphism between
  $V/(V\cap\Sigma)$ and $V\cdot\Sigma/\Sigma$, we see that the
  connected component $W\times T$ of $G/\Sigma$ is isomorphic to
  $V/(V\cap\Sigma)$.  Since the dimension of the maximal torus in
  $V/(V\cap \Sigma)$ is equal to the rank of $V\cap\Sigma$, it follows
  that $\dim T=\rank (V\cap \Sigma)$.  Further, we have
$$S\times F\cong (G/\Sigma)/(G/\Sigma)_0=
(G/\Sigma)/(V\cdot\Sigma/\Sigma)\cong G/V\cdot\Sigma\cong
Z/q(\Sigma),$$ which implies that $\rank S=\rank Z-\rank q(\Sigma)$.
The final assertion follows from the fact that if $G/\Sigma$ is
isomorphic to its dual group, then $\dim T=\rank S$.
\end{proof}

The next lemma closes the main gap in the proof of \cite[Theorem
5.5.13]{ech:mams96}.

\begin{lemma}\label{lemgap}
  Let $G=V\times Z$ be a direct product of a vector group $V$ with a
  finitely generated free abelian group $Z$.  Suppose that
  $(\Sigma_n)_{n\in \NN}$ is a sequence in $\frak K(G)$ which
  converges to $\Sigma\in \frak K(G)$ such that $G/\Sigma_n\cong
  \Sigma_n^{\perp}$ for all $n\in \NN$, and such that $G/\Sigma\cong
  \Sigma^{\perp}$. Then there exists a subsequence of
  $(\Sigma_n)_{n\in \NN}$ (also denoted $(\Sigma_n)_{n\in \NN}$ below)
  which satisfies
  \begin{enumerate}
  \item there exist elements $c_n\in \GL(V)$ such that $c_n\to \id_V$
    in $\GL(V)$ and such that $V\cap \Sigma_n=c_n(V\cap \Sigma)$ for
    all $n\in \NN$.
  \item $q(\Sigma_n)=q(\Sigma)$ for all $n\in\NN$, where $q:G\to Z$
    denotes the quotient map.
  \end{enumerate}
\end{lemma}
\begin{proof} We first claim that $\rank q(\Sigma_n)\geq \rank
  q(\Sigma)$ for almost all $n\in \NN$.  To see this let $s_1,\ldots,
  s_l\in \Sigma$ such that $q(s_1),\ldots, q(s_l)$ is a minimal set of
  generators for $q(\Sigma)$. Since $\Sigma_n\to \Sigma$, it follows
  from the definition of the topology on $\frak K(G)$ that (after
  passing to a subsequence if necessary) there exist sequences
  $(s^i_n)_{n\in \NN}$, $1\leq i\leq l$ such that $s^i_n\in \Sigma_n$
  and $s^i_n\to s_i$ for all $1\leq i\leq l$. But then $q(s^i_n)\to
  q(s_i)$ in $Z$ for all $1\leq i\leq l$, and since $Z$ is discrete we
  have $q(s^i_n)= q(s_i)$ for all $i$ and $n\geq n_0$. Thus
  $q(\Sigma)\subseteq q(\Sigma_n)$ for almost all $n\in \NN$, which
  proves the claim.

  Applying the last assertion of Lemma \ref{lemsplitt} we now get
$$\rank(V\cap \Sigma)=\rank Z-\rank q(\Sigma)\geq\rank Z-\rank
q(\Sigma_n)=\rank(V\cap \Sigma_n).$$ Thus we may apply Lemma
\ref{lemrank} to the sequence $V\cap\Sigma_n$ (which converges to
$V\cap\Sigma$ by Lemma \ref{lemcorrect}) in order to see (after
passing again to a subsequence) that there exists a sequence $c_n\in
\GL(V)$ with $c_n\to \id_V$ and $c_n(V\cap\Sigma)=V\cap \Sigma_n$ for
all $n$.

The only thing which remains to be shown is the assertion that
$q(\Sigma)=q(\Sigma_n)$ for all $n$.  It follows from
$c_n(V\cap\Sigma)=V\cap \Sigma_n$ that
$\rank(V\cap\Sigma_n)=\rank(V\cap\Sigma)$ for all $n\in \NN$, and
therefore
$$\rank q(\Sigma)=\rank Z-\rank(V\cap\Sigma)
=\rank Z-\rank(V\cap\Sigma_n)=\rank q(\Sigma_n)$$ for all $n\in \NN$.
We embed $Z$ as a discrete subgroup of a vector group $W$, so that we
may view all subgroups of $G$ as subgroups of the bigger group
$V\times W$.  Since $\rank\Sigma_n=\rank(V\cap\Sigma_n)+\rank
q(\Sigma_n)$ for all $n$, and the same equation holds for $\Sigma$, we
have $\rank\Sigma_n=\rank \Sigma$ for all $n\in \NN$.  Thus we may
apply Lemma \ref{lemrank} again in order to obtain elements $d_n\in
\GL(V\times W)$ such that $d_n\to \id_{V\times W}$ and
$\Sigma_n=d_n(\Sigma)$ for all $n\in \NN$. Now let $Y\times S$ be a
splitting of $\Sigma$ in its vector group part $Y$ and a finitely
generated free abelian part $S$.  Let ${s_1,\ldots, s_r}$ be a basis
of $Y$ and let ${s_{r+1},\ldots, s_l}$ be a set of generators for $S$.
Then each $\Sigma_n$ has the splitting $\Sigma_n=d_n(Y)\times d_n(S)$
with generating elements $\{d_n(s_1),\ldots, d_n(s_l)\}$.  Let $q$
also denote the projection from $V\times W$ to $W$. Then
$\{q(d_n(s_1)),\ldots, q(d_n(s_l))\}$ is a set of generators of
$q(\Sigma_n)\subseteq Z$ for all $n\in \NN$. Since $d_n\to
\id_{V\times W}$ it follows that $q(d_n(s_i))\to q(s_i)$ for each
$1\leq i\leq l$, which implies, since $Z$ is discrete, that they
eventually coincide. But this proves that $q(\Sigma_n)=q(\Sigma)$ for
all but finitely many $n\in \NN$.
\end{proof}

\begin{proof}[Proof of Theorem \ref{thm}]
  By \cite[Lemma 5.5.6]{ech:mams96} it remains to show that the
  continuity of the symmetrizer map and the type I'ness of the Mackey
  obstructions implies that $h^{\alpha}$ is open as a map onto its
  image.  Each compactly generated Lie group is a quotient of a group
  $V\times Z$ by some discrete subgroup $D$, where $V$ is a vector
  group and $Z$ is a finitely generated free abelian group. Thus, by
  extending the Mackey obstructions to this covering group, we assume
  without loss of generality that $G=V\times Z$.

  Let $(x_n,s_n)_{n\in\NN}$ be a sequence in $\Om\times G$ and
  $(x,s)\in \Om\times G$ such that $h^{\alpha}(x_n,s_n)\to
  h^{\alpha}(x,s)$ in $\Om\times\widehat{G}$. We have to show that,
  after passing to a subsequence if necessary, there exist elements
  $t_n\in\Sigma_{x_n}$ such that $(x_n,t_ns_n)\to (x,s)$ in $\Om\times
  G$.  It is clear that $x_n\to x$ in $\Om$.  Let $q:G\to G/V\cong Z$
  denote the quotient map.  Since $\Sigma_{x_n}\to \Sigma_x$ in $\frak
  K(G)$ and since $s\mapsto h^{\alpha}(x_n,s)$ factors through an
  isomorphism between $G/\Sigma_{x_n}$ and $\Sigma_{x_n}^{\perp}$
  (and, similarly, $G/\Sigma_x\cong \Sigma_x^{\perp}$) by the type I
  assumption on the cocycles, we may apply Lemma \ref{lemgap} in order
  to see that (at least for a subsequence)
  $q(\Sigma_{x_n})=q(\Sigma_x)$ for all $n\in \NN$, and to obtain a
  sequence $c_n\in \GL(V)$ such that $c_n\to \id_V$ and $V\cap
  \Sigma_{x_n}= c_n(V\cap\Sigma_x)$ for all $n\in \NN$.

  For each $n\in \NN$ we now define elements $\sigma_n\in\Aut(G)$ by
  $\sigma_n=c_n\times\id_Z$, where $\id_Z$ denotes the identity on
  $Z$.  Then $\sigma_n\to \id_G$ in $\Aut(G)$ with respect to the
  compact open topology. Moreover, we have
$$\sigma_n(\Sigma_x\cap V)=\Sigma_{x_n}\cap V\quad\text{and}\quad
q(\sigma_n(\Sigma_{x}))=q(\Sigma_{x_n})$$ for all $n\in\NN$.  By
passing to a subsequence, we may assume without loss of generality
that $x_n\neq x_m\neq x$ for all $n,m\in\NN$. Let
$M=\{x_n;n\in\NN\}\cup\{x\}$ and define a map $\tilde{h}:M\times G\to
M\times \widehat{G}$ by
$$\tilde{h}(x_n,s):= (x_n, h_{\om_{x_n}\circ \sigma_n}(s))
\quad\text{and} \quad \tilde{h}(x,s)=h^{\alpha}(x,s)$$ for all
$n\in\NN$ and $s\in G$, where for any $2$-cocycle $\om$ on $G$ and
$\sigma\in \Aut(G)$, $\om\circ\sigma(s,t)=\om(\sigma(s),\sigma(t))$.
Since $h^{\alpha}$ is continuous by \cite[Lemma 5.5.4]{ech:mams96} and
$\sigma_n\to \id_G$ it follows that $\tilde{h}$ is continuous, too,
and we proceed by showing that it is open as a map onto its image.
For this we put $\tilde{h}_{x_n}=h_{\om_{x_n}\circ \sigma_n}$ for all
$n\in\NN$ and $\tilde{h}_x= h_{\om_x}$.  Furthermore, for $n\in\NN$,
let $\tilde{\Sigma}_{x_n}=\sigma_n^{-1}(\Sigma_{x_n})= \ker
\tilde{h}_{x_n}$, and let $\tilde{\Sigma}_x=\Sigma_x=\ker
\tilde{h}_x$.  Then $\tilde{\Sigma}_{x_n}\cap V=\Sigma_x\cap V$ for
all $n\in\NN$, and we may pass over to $G/(\Sigma_x\cap V)$ in order
to show the openness of $\tilde{h}$. After having done this, we may
assume that $\tilde{\Sigma}_{x_n}$ has trivial intersection with the
connected component $G_0$ of $G$ for all $n\in \NN$, and the same may
be assumed for $\tilde{\Sigma}_x$. Moreover, if $q:G\to G/G_0$ denotes
the quotient map, then it follows from our constructions that
$q(\tilde{\Sigma}_{x_n})=q(\tilde{\Sigma}_x)$ for all $n\in\NN$.

Let $H=\tilde{\Sigma}_x\cdot G_0$. Then also $H=\tilde{\Sigma}_y\cdot
G_0$ for all $y\in M$ (since
$q(\tilde{\Sigma}_{x_n})=q(\tilde{\Sigma}_x)$) and we conclude that
$$\tilde{h}_{y}(sH)= 
\tilde{h}_y(s)\tilde{h}_y(\tilde{\Sigma}_y\cdot G_0)=
\tilde{h}_y(s)\tilde{h}_y(G_0)$$ lies in the connected component of
$\widehat{G}$ if and only if $s\in H$.  Suppose now that
$(s_n)_{n\in\NN}$ is a sequence in $G$ such that
$\tilde{h}(x_n,s_n)\to \tilde{h}(x,s)$ for some $s\in G$.  Multiplying
each $s_n$ with $s^{-1}$ we may assume that $\tilde{h}(x_n,s_n)\to
(x,1_G)$ in $M\times \widehat{G}$.  Since $M\times (\widehat G)_0$ is
open in $M\times \widehat G$, it follows from the above observation
that $s_n\in H$ for all but finitely many $n\in \NN$.  Thus, by
multiplying this sequence with $z_n$ in the second variable, for
appropriate elements $z_n\in \tilde{\Sigma}_{x_n}$, we may assume that
$s_n\in G_0$ for all $n\in \NN$ and that $\tilde{h}(x_n,s_n)\to (x,
1_G)$ in $M\times\widehat{G}$.  We want to show that $s_n\to e$ in
$G$.  For this we write $G_0=W\times T$ for a vector group $W$ and a
torus group $T$, and we write $s_n= r_n\cdot t_n$, $r_n\in W$ and
$t_n\in T$.  By the compactness of $T$ we may assume that $t_n\to t$
for some $t\in T$.  It follows that $s_nt_n^{-1}\in W$ for all $n\in
\NN$, and that
$$\tilde{h}_{x_n}(s_nt_n^{-1})=
\tilde{h}_{x_n}(s_n)\tilde{h}_{x_n}(t_n^{-1}) \to
\tilde{h}_{x}(t^{-1}).$$ Since $\tilde{h}_x$ maps the torus part of
$G$ into the torus part of $\widehat{G}$, it follows that
$\tilde{h}_{x}(t^{-1})|_W=1_W\in\widehat{W}$.  Thus
$\tilde{h}_{x_n}(s_nt_n^{-1})|_W\to 1_W$ in $\widehat{W}$.  Since each
$\tilde{\Sigma}_{x_n}$ has trivial intersection with $G_0$, and hence
also with $W$, we see that the restriction of $\om_{x_n}\circ\sigma_n$
to $W\times W$ is totally skew for all $n\in\NN$, and the same applies
to the restriction of $\om_x$ to $W\times W$.  Thus, for each $y\in
M$, we obtain an isomorphism $W\to\widehat{W}$ given by $s\mapsto
\tilde{h}_y(s)|_W$.  By identifying $W$ with $\widehat{W}$ via any
fixed linear isomorphism, the continuity of $\tilde{h}$ implies that
the map $M\to \GL(W); y\mapsto \tilde{h}_y|_W$ is continuous, too.
Thus, the map $M\to \GL(W); y\mapsto (\tilde{h}_y|_W)^{-1}$ is also
continuous.  This implies that $s_nt_n^{-1}\to \{0\}$ in $W$, which in
turn yields $s_n\to t$ in $G$. But, since by assumption
$\tilde{h}_{x_n}(s_n)\to \tilde{h}_x(t)=1_G$ in $\widehat{G}$, and
since $\ker \tilde{h}_x\cap G_0=\Sigma_x\cap G_0=\{e\}$, it follows
that $t=1$. This completes the proof for the openness of $\tilde{h}$.

Finally, let $(s_n)_{n\in\NN}$ be a sequence in $G$ such that
$h^{\alpha}(x_n,s_n)\to h^{\alpha}(x,s)$ for some $s\in G$.  Then it
follows that $\tilde{h}(x_n, \sigma_n^{-1}(s_n))=
(x_n,h_{\om_{x_n}}(s_n)\circ\sigma_n)$ converges to
$(x,h_{\om_x}(s))=\tilde{h}(x,s)$ in $M\times\widehat{G}$.  Since
$\tilde{h}$ is open as a map onto its image, we may multiply each
$\sigma_n(s_n)$ with an appropriate element $r_n\in
\tilde{\Sigma}_{x_n}=\sigma_n^{-1}(\Sigma_{x_n})$ such that, after
possibly passing to a subsequence, $(x_n, \sigma_n^{-1}(s_n)r_n)\to
(x,s)$ in $M\times G$.  Thus, if $r_n'=\sigma_n(r_n)$, we have
$r_n'\in\Sigma_{x_n}$ for all $n\in\NN$ and $(x_n,s_nr_n')\to (x,s)$
in $M\times G\subseteq \Om\times G$.  This completes the proof.
\end{proof}

\bibliographystyle{amsxport} 

\begin{bibdiv}
\begin{biblist}

\bib{arv:invitation}{book}{
      author={Arveson, William},
       title={An {I}nvitation to {$C\sp*$}-algebras},
   publisher={Springer-Verlag},
     address={New York},
        date={1976},
        note={Graduate Texts in Mathematics, No. 39},
      review={\MR{MR0512360 (58 \#23621)}},
}

\bib{bprs:nyjm00}{article}{
      author={Bates, Teresa},
      author={Pask, David},
      author={Raeburn, Iain},
      author={Szyma{\'n}ski, Wojciech},
       title={The {$C\sp *$}-algebras of row-finite graphs},
        date={2000},
        ISSN={1076-9803},
     journal={New York J. Math.},
      volume={6},
       pages={307\ndash 324 (electronic)},
      review={\MR{MR1777234 (2001k:46084)}},
}

\bib{deicke}{article}{
	 author={Deicke, Klaus},
       title={Pointwise unitary coactions on $C\sp *$-algebras with continuous trace},
        date={2000},
        ISSN={0379-4024},
     journal={J. Operator Theory},
      volume={43},
       pages={295\ndash 327},
      review={\MR{MR1753413 (2001d:46098)}},
}


\bib{dix:von}{book}{
      author={Dixmier, Jacques},
       title={von {N}eumann algebras},
      series={North-Holland Mathematical Library},
   publisher={North-Holland Publishing Co.},
     address={Amsterdam},
        date={1981},
      volume={27},
        ISBN={0-444-86308-7},
        note={With a preface by E. C. Lance, Translated from the second French
  edition by F. Jellett},
      review={\MR{MR641217 (83a:46004)}},
}

\bib{ech:mams96}{article}{
      author={Echterhoff, Siegfried},
       title={Crossed products with continuous trace},
        date={1996},
        ISSN={0065-9266},
     journal={Mem. Amer. Math. Soc.},
      volume={123},
      number={586},
       pages={viii+134},
      review={\MR{98f:46055}},
}

\bib{eff:tams63}{article}{
      author={Effros, Edward~G.},
       title={A decomposition theory for representations of {$C^*$}-algebras},
        date={1963},
     journal={Trans. Amer. Math. Soc.},
      volume={107},
       pages={83\ndash 106},
      review={\MR{26 \#4202}},
}

\bib{fol:course}{book}{
      author={Folland, Gerald~B.},
       title={A course in abstract harmonic analysis},
      series={Studies in Advanced Mathematics},
   publisher={CRC Press},
     address={Boca Raton, FL},
        date={1995},
        ISBN={0-8493-8490-7},
      review={\MR{MR1397028 (98c:43001)}},
}

\bib{gli:pjm62}{article}{
      author={Glimm, James},
       title={Families of induced representations},
        date={1962},
        ISSN={0030-8730},
     journal={Pacific J. Math.},
      volume={12},
       pages={885\ndash 911},
      review={\MR{MR0146297 (26 \#3819)}},
}

\bib{gooros:im79}{article}{
      author={Gootman, Elliot~C.},
      author={Rosenberg, Jonathan},
       title={The structure of crossed product {$C^*$}-algebras: a proof of the
  generalized {E}ffros-{H}ahn conjecture},
        date={1979},
        ISSN={0020-9910},
     journal={Invent. Math.},
      volume={52},
      number={3},
       pages={283\ndash 298},
      review={\MR{80h:46091}},
}

\bib{gre:am78}{article}{
      author={Green, Philip},
       title={The local structure of twisted covariance algebras},
        date={1978},
     journal={Acta Math.},
      volume={140},
       pages={191\ndash 250},
}

\bib{kprr:jfa97}{article}{
      author={Kumjian, Alex},
      author={Pask, David},
      author={Raeburn, Iain},
      author={Renault, Jean},
       title={Graphs, groupoids, and {C}untz-{K}rieger algebras},
        date={1997},
        ISSN={0022-1236},
     journal={J. Funct. Anal.},
      volume={144},
      number={2},
       pages={505\ndash 541},
      review={\MR{MR1432596 (98g:46083)}},
}

\bib{kumpas:etds99}{article}{
      author={Kumjian, Alex},
      author={Pask, David},
       title={{$C\sp *$}-algebras of directed graphs and group actions},
        date={1999},
        ISSN={0143-3857},
     journal={Ergodic Theory Dynam. Systems},
      volume={19},
      number={6},
       pages={1503\ndash 1519},
      review={\MR{MR1738948 (2000m:46125)}},
}

\bib{mac:pnasus49}{article}{
      author={Mackey, George~W.},
       title={Imprimitivity for representations of locally compact groups.
  {I}},
        date={1949},
     journal={Proc. Nat. Acad. Sci. U. S. A.},
      volume={35},
       pages={537\ndash 545},
      review={\MR{MR0031489 (11,158b)}},
}

\bib{olerae:jfa90}{article}{
      author={Olesen, Dorte},
      author={Raeburn, Iain},
       title={Pointwise unitary automorphism groups},
        date={1990},
        ISSN={0022-1236},
     journal={J. Funct. Anal.},
      volume={93},
      number={2},
       pages={278\ndash 309},
      review={\MR{92b:46105}},
}

\bib{rw:morita}{book}{
      author={Raeburn, Iain},
      author={Williams, Dana~P.},
       title={Morita equivalence and continuous-trace {$C^*$}-algebras},
      series={Mathematical Surveys and Monographs},
   publisher={American Mathematical Society},
     address={Providence, RI},
        date={1998},
      volume={60},
        ISBN={0-8218-0860-5},
      review={\MR{2000c:46108}},
}

\bib{sau:ma77}{article}{
      author={Sauvageot, Jean-Luc},
       title={Id\'eaux primitifs de certains produits crois\'es},
        date={1977/78},
        ISSN={0025-5831},
     journal={Math. Ann.},
      volume={231},
      number={1},
       pages={61\ndash 76},
      review={\MR{MR473355 (80d:46112)}},
}

\bib{sau:jfa79}{article}{
      author={Sauvageot, Jean-Luc},
       title={Id\'eaux primitifs induits dans les produits crois\'es},
        date={1979},
        ISSN={0022-1236},
     journal={J. Funct. Anal.},
      volume={32},
      number={3},
       pages={381\ndash 392},
      review={\MR{81a:46080}},
}

\bib{wea:jfa03}{article}{
      author={Weaver, Nik},
       title={A prime {$C\sp *$}-algebra that is not primitive},
        date={2003},
        ISSN={0022-1236},
     journal={J. Funct. Anal.},
      volume={203},
      number={2},
       pages={356\ndash 361},
      review={\MR{MR2003352 (2004g:46075)}},
}

\bib{wil:tams81}{article}{
      author={Williams, Dana~P.},
       title={The topology on the primitive ideal space of transformation group
  {$C\sp{\ast} $}-algebras and {C}.{C}.{R}. transformation group {$C\sp{\ast}
  $}-algebras},
        date={1981},
        ISSN={0002-9947},
     journal={Trans. Amer. Math. Soc.},
      volume={266},
      number={2},
       pages={335\ndash 359},
      review={\MR{MR617538 (82h:46081)}},
}

\bib{zel:jmpa68}{article}{
      author={Zeller-Meier, G.},
       title={Produits crois\'es d'une {$C\sp{\ast} $}-alg\`ebre par un groupe
  d'automorphismes},
        date={1968},
        ISSN={0021-7824},
     journal={J. Math. Pures Appl. (9)},
      volume={47},
       pages={101\ndash 239},
      review={\MR{MR0241994 (39 \#3329)}},
}

\end{biblist}
\end{bibdiv}

\def\noopsort#1{}\def\cprime{$'$} \def\sp{^}

\end{document}